\documentclass[times,10pt]{article}

\begin{document}

\newtheorem{Remark}{Remark}[section]
\newtheorem{Lemma}{Lemma}[section]
\newtheorem{Proposition}{Proposition}[section]
\newtheorem{Theorem}{Theorem}[section]
\newtheorem{Exercise}{Exercise}[section] 
\newtheorem{Corollary}{Corollary}[section] 
\newtheorem{Definition}{Definition}[section]
\newtheorem{Example}{Example}[section]
\def\ad{\mbox{ad\,}} \def\tr{\mbox{tr\,}}\def\char{\mbox{char\,}}
\def\mod{\mbox{mod\,}}\def\Ad{\mbox{Ad\,}}
\def\diag{\mbox{diag\,}}\def\ad{\mbox{ad\,}} \def\tr{\mbox{tr\,}} 
\def\End{\mbox{End\,}}\def\GL{\mbox{GL\,}}\def\SL{\mbox{SL\,}}
\def\Out{\mbox{Out\,}}\def\Int{\mbox{Int\,}}
\def\Mat{\mbox{Mat\,}}\def\Hom{\mbox{Hom\,}}\def\Iso{\mbox{Iso\,}}
\def\Aut{\mbox{Aut\,}}\def\gl{\mbox{gl\,}}\def\sl{\mbox{sl\,}}
\def\o{\mbox{o\,}}\def\sp{\mbox{sp\,}}
\def\im{\mbox{im\,}}\def\ker{\mbox{ker\,}}\def\deg{\mbox{deg\,}}
\def\id{\mbox{id\,}}\def\mod{\mbox{mod\,}}
\def\grad{\mbox{\rm grad\,}}\def\rot{\mbox{\rm rot\,}}
\def\div{\mbox{\rm div\,}}\def\Grad{\mbox{\rm Grad\,}}
\def\det{\mbox{\rm det\,}} \def\ctg{\mbox{\rm ctg\,}}\def\tg{\mbox{\rm tg\,}}
\def\sn{\mbox{\rm sn\,}}\def\th{\mbox{\rm th\,}}\def\dn{\mbox{\rm dn\,}}
\def\pd#1,#2{\frac{\partial#1}{\partial#2}}\def\diag{\mbox{\rm diag\,}}
\def\v#1{\overline{#1}}\def\sh{\mbox{\rm sh\,}}\def\ch{\mbox{\rm ch\,}}
\def\d#1,#2{\frac{d#1}{d#2}}\def\th{\mbox{\rm th\,}}
\def\qed{\hbox{${\vcenter{\vbox{
\hrule height 0.4pt\hbox{\vrule width 0.4pt height 6pt
\kern 5pt\vrule width 0.4pt}\hrule height 0.4pt}}}$}}

\begin{titlepage} 
\title{Lie algebra extensions related with linear bundles of Lie
brackets} 
\author{A  B  Yanovski \thanks{\it On leave of absence from the Faculty of
Mathematics and Informatics, St.  Kliment Ohridski University, James
Boucher Blvd 1164 - Sofia, Bulgaria}}

\maketitle 
\begin{center} 
Universidade Federal de Sergipe, Cidade Universit\'aria,\\ ''Jos\'e
Alo\'{\i}sio de Campos'', 49.100-000-S\~ao Christ\'ov\~ao, SE-Brazil
\end{center} 
\begin{abstract} 
We consider some special type extensions of an arbitrary Lie algebra
${\cal G}$, arising in the theory of Lie-Poisson structures over
$({\cal G}^*)^n$, where ${\cal G}^*$ is the dual of ${\cal G}$. We
show that some classes of these extensions can be constructed in a
natural way using some linear bundles of Lie algebras.
\end{abstract}
\begin{center} 
37K30,  37K05,  17B80
\end{center} 

\end{titlepage}
\section{Introduction}
Let ${\cal G}$ be an arbitrary Lie algebra over one of the classical
fields, which we shall denote by ${\bf K}$ and let us consider the
linear space ${\cal G}^n$ with elements 
\begin{equation}
{\bf x}=(x_1,x_2,\ldots,x_n),\quad x_i\in {\cal G}.
\end{equation}
In this paper, following \cite{JLT}, we shall be interested in an
universal manner of defining Lie brackets over ${\cal G}^n$. More
specifically, we shall consider brackets with the following properties:
\begin{itemize}
\item The brackets have the form:
\begin{equation}\label{eq:DEF}
([{\bf x},{\bf y}]_W)_s=\sum\limits_{i,j=1}^n W^{ij}_s[x_i,y_j];
\quad W^{ij}_s\in {\bf K}.
\end{equation}
\item For $W=(W^{ij}_s)$ fixed, (\ref{eq:DEF}) is a Lie bracket for
arbitrary Lie algebra ${\cal G}$.
\end{itemize}

It is not difficult to see that $[{\bf x},{\bf y}]_W$ is a Lie
bracket of the type desired if the tensor $W=(W^{ij}_s)$
satisfy the conditions:
\begin{equation}\label{eq:0}
\begin{array}{c}
W^{ij}_s=W^{ji}_s\\[4pt]
\sum\limits_{k=1}^n(W_i^{sk}W_{k}^{qp}-W_i^{qk}W_{k}^{sp})=0.
\end{array}
\end{equation}

Algebraic structures as described above exist, that is there exist solutions 
to (\ref{eq:0}). For example if $W^{ij}_s=\delta_{s}^i\delta_{s}^j$ we get 
$[{\bf  x},{\bf y}]_s=[x_s,y_s]$ which corresponds to the direct sum 
structure:
\begin{equation}
{\cal G}^n=\oplus_{i=1}^n {\cal G}.
\end{equation}
If not assumed something else we shall always imply that ${\cal G}^n$
is a Lie algebra with the direct sum structure. However, another, and
quite interesting examples of Lie algebra structures defined by
tensors of the type $W_{s}^{ij}$ exist and already have been
applied to describe the Poisson structures of some hydrodynamic and
magneto-hydrodynamic models, see \cite{JLT}, where there is extensive
bibliography and \cite{JLT1} in order to describe a Poisson structure
in one interesting mechanical model. We imply here the Poisson structures
on the coalgebras $({\cal G}^n)^*$ - the so called Poisson-Lie, or
Kirillov structures, see \cite{Lich,Ki,MR}. It is known that they
play important role in the Hamiltonian formalism and in the modern
theories of integration of evolution equations,\cite{TrFo, ReyST2}. 
As the Lie algebra structure on ${\cal G}^n$ is not unique we shall
use the symbol ${\cal G}^n_W$ for the structure defined by the tensor
$W_k^{ij}$, or shall say explicitly what we mean.

In the case ${\bf K}={\bf C}$ the tensors $W_{s}^{ij}$ and the
corresponding Lie algebra structures have been given a comprehensive
analysis, see \cite{JLT}. For convenience of the reader we shall
remind some facts of this analysis below. First of all however, for
reasons that will become clear later, we scale the indices in a different way and we
assume that they run from $0$ to $n$ instead from $1$ to $n$.

The starting point of the analysis is the fact that if we introduce
the $n+1$ matrices $W^{(k)}$, $0\leq k\leq n$ with components
$(W^{(k)})_i^j=W^{kj}_i$ then the second equation in (\ref{eq:0})
means that the matrices $W^{(k)}$ commute. It follows by a linear
transformation  that all ${W}^{(s)}$ can be put simultaneously into a
block-diagonal form, each block being low-triangular with the
generalized eigenvalues on the diagonal. The block structure
corresponds to a splitting of the algebra ${\cal G}_W^n$ into a
direct sum.  Thus if we limit ourselves with the irreducible case, we
can assume that all the matrices $W^{(s)}$ are low triangular.  Now,
the symmetry of $W_{i}^{jk}$ entails that the generalized eigenvalues
of $W^{(s)}$ for $s>0$ are zero and so for $s>0$ the matrices
$W^{(s)}$ are low-triangular with zeroes on the diagonal, and hence
are nilpotent.  For the generalized eigenvalue $\mu_0$ of $W^{(0)}$
there are two alternatives: $\mu_0\neq 0$ and $\mu_0=0$. If
$\mu_0\neq 0$ it can be shown that with a suitable linear transform
we can achieve $W^{(0)}$ to be equal to the unit matrix conserving the
lower-triangular form of the rest of $W^{(i)}$. This case is called
in \cite{JLT} semisimple. When $\mu_0=0$ one can see that $W^{(n)}$
is the zero matrix and this case is called solvable case. There is
correspondence between solvable and semisimple extensions and one can
limit ourself with the study of the solvable ones. The correspondence
is very simple, having the semisimple case with tensor $W_i^{jk}$,
$0\leq i,j,k\leq n$ and the matrices $W^{(i)}$, $0\leq i\leq n$ we
just consider the matrices $\bar{W}^{(i)}$, $1\leq i\leq n$ where
$\bar{W}^{(i})$ is obtained from $W^{(i})$ cutting the first row and
column, or in other words we consider $\bar{W}^{ij}_k=W^{ij}_k$ for
$1\leq i,j,k\leq n$. The operation will be referred as ''finding the
semisimple part''. We shall always assume that in the semisimple case
indices run from $0$ to some $m$ and in the solvable case from $1$ to
some $k$.

If ${\cal G}_W^n$ is solvable extension, the vector space ${\cal
G}^n$ can be split into
\begin{equation}
{\cal G}^n={\cal F}_n^{(1)}\oplus{\cal F}_n^{(1)}\oplus\ldots
\oplus {\cal F}_n^{(n)}
\end{equation}
where ${\cal F}_n^{(i)}$ consists of those ${\bf x}$ for which
$x_j=0$ for all $j\neq i$. For $1\leq k\leq n$ let us set 
\begin{equation}
{\cal F}[n,k]={\cal F}_n^{(k)}\oplus {\cal F}_n^{(k+1)}\oplus\ldots 
\oplus {\cal F}_n^{(n)}
\end{equation}
and for $k>n$ let ${\cal F}[n,k]=0$. Obviously,
\begin{equation}
0\subset {\cal F}^{(n)}_n={\cal F}[n,n]\subset {\cal F}[n,n-1]\subset
\ldots\subset {\cal F}[n,1]={\cal G}_W^n
\end{equation}
It can be shown, that in the solvable case the coordinates can be
chosen in such way, that
\begin{equation}
\left[{\cal F}[n,k],{\cal F}[n,s]\right]\subset {\cal F}[n,(\max(k,s)) +1]
\end{equation}
In particular ${\cal F}_n^{(n)}$ is Abelian ideal. Denote by ${\cal
S}^{k}_{n}$ the maximal Abelian ideal of the type ${\cal
F}[n,n-k+1]$. As ${\cal F}[n,n]$ is an Abelian ideal ${\cal
S}^{n}_{k}$ is not $0$. Then we have the exact sequence
\begin{equation}
0\mapsto {\cal S}_{n}^{k}\mapsto {\cal G}^n_W\mapsto {\cal G}^n_W\slash
{\cal S}^{n}_{k}\mapsto 0
\end{equation}
which defines an extension through an Abelian kernel, see \cite{GoGr,ChEi},
in the sense that it is usually understood, which explains why the
structures are called extensions.

In the present article we shall give a very simple
construction of a new class of $W$ tensors and show how one can
use them to construct in natural way the so called Leibnitz
extensions which seem to be the only regular class of universal
extensions, that is they appear for every $n$. They are more
interesting, because as it is shown in \cite{JLT} they are in some
sense maximal and are given by the very simple formula
$W^{ij}_k=\delta_k^{i+j}$ where $\delta_s^m$ is the Kronecker symbol.
Thus we believe we can shade some light on the origin of the
algebraic structures described in the above. To this end and for
convenience of the reader we shall recall some facts and
constructions from the theory of the so-called linear bundles of Lie
algebras (LBLA), see \cite{Ya3}, which we introduce below.

\section{Linear bundles of Lie algebras}
Let us start with some definitions.
\begin{Definition}
Let ${\cal G}$ and $V$ be  vector spaces and let for arbitrary 
$v\in V$ is defined the Lie bracket on ${\cal G}$:
\begin{equation}
(X,Y)\rightarrow L_v(X,Y)\equiv [X,Y]_v\equiv \ad_X^v(Y)
\end{equation}
such that for $v_1,v_2\in V,\quad a_1,a_2\in {\bf K}$
\begin{equation}
[X,Y]_{a_1v_1}+[X,Y]_{a_2v_2}=[X,Y]_{a_1v_1+a_2v_2}.
\end{equation}
We say that $({\cal G},V)$ is a linear bundle of Lie algebras and
if $V$ is finite dimensional we shall call the dimension of $V$  dimension 
of the linear bundle  $({\cal G},V)$.

The bundle $({\cal G},{\bf K})$ defined as:
\begin{equation}
[X,Y]_a=a[X,Y],\quad X,Y\in {\cal G},~a\in {\bf K}
\end{equation}
is called the trivial bundle corresponding to ${\cal G}$.
\end{Definition}
The space ${\cal G}$ endowed with the bracket $L_v$ (or $[~,~]_v$)
shall be denoted by ${\cal G}_v$. When ${\cal G}_v$ coincides with
some Lie algebra with respect to the usual commutator we don't write
any index. We denote by the same letter the algebra ${\cal G}$ and
the underlying vector space. Denote by $d^{\rho}$ the
coboundary operator for ${\cal G}$ with respect to the
representation $({\cal G},\rho, W)$ ($W$-vector space,
$\rho$-representation of ${\cal G}$ in $W$) and the corresponding
spaces of cocycles, coboundaries and cohomologies by 
$Z_{\rho}^{k}({\cal G}_v,W)$,
$B_{\rho}^{k}({\cal G}_v,W)$, $H^{k}_{\rho}({\cal G}_v,W)$.  If the adjoint
representation of ${\cal G}$ or the trivial representation is
implied we don't write the symbol $\rho$. The complexes are of course
the graded modules of skewsymmetric maps from ${\cal G}$ into $W$ and
the cohomologies are defined as usual, see \cite{ChEi,GoGr}. When we
deal with algebras ${\cal G}_v$ of some bundle we label the
corresponding cohomology groups  and the coboundary operator by $v$.
Thus we write $d_v$ for the
coboundary operator defined by the adjoint action of ${\cal G}_v$.
It is easy to check that the fact that for $u,v\in V$ the expression
$L_u+L_v$ is Lie bracket entails that for $X_1,X_2,X_3\in {\cal G}$
one has the identity:
\begin{equation}\label{eq:osn}
[X_1,[X_2,X_3]_u]_v+[X_1,[X_2,X_3]_v]_u+cycl(1,2,3)=0.
\end{equation}
The notation  $cycl(1,2,3)$ means that one must add to the first two
terms expressions obtained from them by cyclic permutation of the
indices 1,2,3. 
Taking into account (\ref{eq:osn}) one can see that if $V$ is a
vector space and $L_v$, $v\in V$ is a family of Lie brackets having
the property that $L_{\mu v}=\mu L_v$ it will be linear bundle of Lie
brackets if and only if
\begin{equation}\label{eq:h}
d_vL_w=d_wL_v=0,\qquad v,w\in V.
\end{equation}
This means that
\begin{Proposition}\label{Theorem:TF}
For arbitrary $v,w\in V$ the bracket $L_w$ is $2$-cocycle for the
coboundary operator $d_v$. 
\end{Proposition}
More generally
\begin{Proposition}\label{Proposition:I}
If $[X,Y]_1$ and $[X,Y]_2$ are different Lie brackets on ${\cal G}$ then their
sum $[X,Y]_1+[X,Y]_2$ is a Lie bracket if and only if (\ref{eq:h})
holds (with indices $1,2$ instead of $v,w$). In this case for arbitrary
$\lambda,\mu\in {\bf K}$ the expression:
\begin{equation}\label{eq:fml}
\lambda [X,Y]_1+\mu [X,Y]_2
\end{equation}
is a Lie bracket.
\end{Proposition}
In case the Lie brackets $L_1,L_2$ and $L_1+L_2$ are Lie brackets
they are called compatible and it is clear, that they define
$2$-dimensional linear bundle of Lie algebras over ${\cal G}$.

The applications of the Lie bundles of Lie algebras to Dynamics of
Hamiltonian systems are consequence
of the following simple proposition, which applied to the corresponding
Poisson-Lie  tensors allows to construct sets of functions in
involution, see \cite{TrFo}.
\begin{Proposition}\label{Proposition:ZZ}
Denote by $Z({\cal G}_w)$ the center of ${\cal G}_w$, $w\in V$. We have:
\begin{itemize}
\item $Z({\cal G}_u)$ is subalgebra in all the algebras ${\cal
G}_v,~v\in V$. 
\item If $X_1\in Z({\cal G}_u)$ and $X_2\in Z({\cal G}_v)$ then
$[X_1,X_2]_u, [X_1,X_2]_v\in Z({\cal G}_w)$ for $w=\lambda
u+\mu v$;  $\lambda,\mu$ being arbitrary numbers.
\end{itemize}
\end{Proposition}
Our intention however is another, we shall try to find solutions
of the equations (\ref{eq:DEF}) using some of the typical constructions
of LBLA theory in order to find new structures and deforming these
structures into another ones.

\section{Construction of Lie algebra structures over ${\cal G}^n$}

Let ${\cal G}$ be finite dimensional Lie algebra. Let ${\cal A}$ be an 
associative algebra over ${\bf K}$ with unity ${\bf 1}$
such that ${\cal G}$ is imbedded in a natural way in ${\cal A}$, that is
$[x,y]=xy-yx$. (We denote by the same letters the elements of ${\cal
G}$ and their images in ${\cal A}$).  For example, the role of ${\cal
A}$ can be played by the universal enveloping algebra $U({\cal G})$ of ${\cal
G}$. Also, if ${\cal G}\subset \Mat(N,{\bf K})$ is a matrix algebra
we can set ${\cal A}=\Mat(N,{\bf K})$ - the associative algebra of
$N\times N$ matrices.

Consider now $n\times n$ matrices with entries in ${\cal A}$. The
set of these matrices form an associative algebra with
unity which we shall denote by $\Mat(n,{\cal A})$. For example, when
${\cal A}=\Mat(n,{\bf K})$ the algebra $\Mat(n,{\cal A})$ is isomorphic to
$\Mat(nN,{\bf K})$. Let us take the following imbedding of ${\cal
A}^n$ (as a vector space) in $\Mat(n,{\cal A})$:
\begin{equation}\label{eq:1}
{\cal A}^n\ni {\bf x}\mapsto X\in \Mat(n,{\cal A}),\quad
X_{ij}=x_{(i+j) mod(n)}.
\end{equation}
Of course, we then have also an imbedding of ${\cal G}^n$ in
$\Mat(n,{\cal A})$ too.

It is relevant to assume that all the indices are $\mod (n)$ and
hence they run over $0,1,2,\ldots n-1$. The matrix $X$
introduced in (\ref{eq:1}) has the form:
\begin{equation}\label{eq:2}
X=\left(\begin{array}{cccccc}
x_0&x_1&x_2&\ldots&x_{n-2}&x_{n-1}\\
x_1&x_2&x_3&\ldots&x_{n-1}&x_0\\
\ldots&\ldots&\ldots&\ldots&\ldots&\ldots\\
x_{n-1}&x_0&x_1&\ldots&x_{n-3}&x_{n-2}
\end{array}\right), \qquad x_i\in {\cal A}.
\end{equation}
The set of matrices having the same structure with $x_i\in {\cal
A}$ will be denoted by ${\cal A}^n$. Analogously
by ${\cal G}^n$  we shall denote the set of matrices with the form
(\ref{eq:2}) with $x_i\in {\cal G}$.  Now we start by defining some Lie
natural algebra structures over ${\cal A}^n$ and then restrict them
to ${\cal G}^n$.
\begin{Proposition}
If $X,Y,A\in {\cal A}^n$, then $XAY$ (the polytop product) 
belongs to ${\cal A}^n$ and endows ${\cal A}^n$ with a structure of
an associative algebra.
\end{Proposition}
\begin{Corollary}
For fixed $A\in {\cal A}^n$ the bracket
\begin{equation}
[X,Y]_A=XAY-YAX
\end{equation}
defines on ${\cal A}^n$ a structure of a Lie algebra.
\end{Corollary}
Thus $({\cal A}^n,{\cal A}^n)$ is a typical LBLA, see \cite{Ya3}.
We can easily calculate the element $[{\bf x},{\bf y}]_A\in {\cal
A}^n$ that corresponds to $[X,Y]_A$:
\begin{equation}\label{eq:3}
\begin{array}{c}
([{\bf x},{\bf y}]_A)_i=\mod(n)\sum\limits_{s,k=0}^{n-1}(x_sa_{s+k-i}y_k-
y_ka_{s+k-i}x_s)=\\[4pt]
\mod(n)\sum\limits_{s,k=0}^{n-1}[x_s,y_k]_{a_{s+k-i}},\quad
[x,y]_a\equiv xay-yax,
\end{array}
\end{equation}
where the symbol $\mod(n)$ before the sum remind us that 
in the final expressions all the indices are understood modulo $n$ and thus take values from
$0$ to $n-1$. 

It is not easy to give a general receipt how to restrict the above
structures onto ${\cal G}^n$, but if $({\cal G},V)$ is  a linear
bundle of Lie algebras with the bracket $[x,y]_a=xay-yax$ and
$V\subset {\cal A}$ then for $A\in V^n$ the formula (\ref{eq:3}) will
give a linear bundle of Lie algebras, which we shall denote by
$({\cal G}^n,V^n)$. For example, we can take ${\cal G}=so(p,{\bf K})$
- the algebra of the skew-symmetric matrices and ${\cal
A}=\Mat(p,{\bf K})$. For $x_i\in
so(p,{\bf K})$ $X$ will be skew-symmetric $pn\times pn$ matrix over
${\bf K}$. We can take $A$ to have the same structure as that of $X$
with $a_i$ symmetric matrices ($a_i\in sym(n,{\bf K})$) and then $A$
will be symmetric $pn\times pn$ matrix over ${\bf K}$. As
$(so(p,{\bf K}),sym(p,{\bf K}))$ is LBLA we shall obtain in this
manner the LBLA $(so(p,{\bf K})^n,sym(p,{\bf K})^n)$.

But one can easily see that the construction works even if $({\cal
G},V)$ is the trivial bundle, that is $V={\bf K}$. In this case the
bracket in the bundle $({\cal G},{\bf K})$ is simply
$[x,y]_a=a[x,y]$, $a\in {\bf K}$. We then assume that the components
$a_i$ of $A$ are of the type $a_i=\alpha_i{\bf 1}$,
$\alpha_i\in {\bf K}$. Then we have the bracket
\begin{equation}\label{eq:4}
([{\bf x},{\bf y}]^{\bf \alpha})_i=\mod(n)\sum\limits_{s,k=0}^{n-1}
\alpha_{s+k-i}[x_s,y_k],\qquad ({\bf\alpha})_i=\alpha_i
\end{equation}
and it will be exactly of the type we are looking for. Thus we obtain a LBLA
$({\cal G}^n,{\bf K}^n)$ and have found a class of solutions of (\ref{eq:0}),
given by the formula:
\begin{equation}
W[\alpha]^{sk}_i=\alpha_{(s+k-i)mod(n)}.
\end{equation} 
To the best of our knowledge, though the above solutions are very
simple, they have not been considered.

We want to make now an important remark. We have used
the fact that the Lie algebra ${\cal G}$ is finite
dimensional only in order to ensure the existence of ${\cal A}$.  But
now, after $W[\alpha]_i^{sk}$ is found we do not need any more this
requirement and can use the tensor $W[\alpha]_i^{sk}$ for arbitrary
${\cal G}$.

Let us consider with more detail the structures defined in the above.
In what follows we denote by $e_i\in {\bf K}^n$ the vector with
components $(e_i)_j=\delta_i^j$ and by $W[\alpha]^{(k)}$ the matrix
with components $(W[\alpha]^{(k)})_i^j=W[\alpha]^{kj}_i$ with respect
to the basis $\{e_i\}_0^{n-1}$. In order to simplify the notations,
when $\alpha=e_0$ we shall omit the symbol $e_0$ and shall write
simply $W^{jk}_i$ and $W^{(k)}$.  

It is not difficult to check that
\begin{equation}
W[\alpha]^{(s)}=\sum\limits_{k=0}^{n-1}\alpha_kW^{(s-k)}
\end{equation}
(all indices are $\mod(n)$).

For the reasons that will become obvious later we shall assume that ${\bf
K}={\bf C}$. Next, it is not difficult to show that the matrices $W^{(p)}$,
$p=0,1,\ldots, n-1$ in act on the basis $\{e_i\}_0^{n-1}$ as follows:
\begin{equation}
W^{(p)}e_i=e_{(i+p)mod(n)}.
\end{equation}
Therefore $W^{(0)}={\bf 1},~W^{(p)}=(W^{(1)})^p$ and see that
the matrices $W^{(p)}[\alpha]$, $p=0,1,\ldots, n-1$ commute not only for
fixed $\alpha$ (which is obvious from the properties of the tensors
$W[\alpha]^{jk}_i$) but for all the choices of $\alpha$. It is also
not difficult to see that we can simultaneously diagonalize all
$W[\alpha]^{(p)}$.  For this we take into account that
$W^{(p)}f_l=\omega^{lp}f_l$ where $\omega=e^{2\pi i\over n}$ and
\begin{equation}
f_l={1\over n}\sum\limits_{j=0}^{n-1}\omega^{-jl}e_j,\qquad l=0,1,\ldots,n-1.
\end{equation}
Thus the diagonalization of $W^{(p)}$ is performed by similarity transform
\begin{equation}
W^{(p)}\mapsto \Omega^{-1}W^{(p)}\Omega=
\diag(1,\omega^p,\omega^{2p},\ldots \omega^{(n-1)p}),
\end{equation}
where $\Omega_i^j={1\over n}\omega^{-ij}$.

The tensor $W_i^{jk}$ then transforms as follows:
\begin{equation}
\begin{array}{c}
(W_k^{ij})\mapsto (\tilde{W}_k^{ij})\\[4pt]
\tilde{W}_k^{ij}=\sum\limits_{s,p,q=0}^{n-1}
W_s^{pq}\Omega^p_i\Omega^q_j(\Omega^{-1})_k^s=
\sum\limits_{p=0}^{n-1}
(\Omega^{-1}W^{(p)}\Omega)_k^j\Omega^i_p=\delta_k^i\delta_k^{j}.
\end{array}
\end{equation}
As we know $\tilde{W}_k^{ij}=\delta_k^i\delta_k^{j}$ corresponds to
the direct sum structure.
The case $W[\alpha]^{jk}_i$ is quite analogous, with similar
calculations we get
\begin{equation}\label{eq:mus}
\begin{array}{c}
\tilde{W}[\alpha]^{jk}_i=\mu_i\tilde{W}^{jk}_i,\\[4pt]
\mu_i=\sum\limits_{r=0}^{n-1}\alpha_r\omega^{-ri}.
\end{array}
\end{equation}
Thus, the resulting structure is equivalent to:
\begin{equation}
[\tilde{\bf x},\tilde{\bf y}]_i=\mu_i[\tilde{x}_i,\tilde{y}_i].
\end{equation}
In other words, we have:
\begin{Proposition}
Let the numbers $\mu_i,~0\leq i\leq n-1$ are defined as in
(\ref{eq:mus}) and let $n-m$ be the total number of $\mu_i$ that are
equal to zero.  Then the Lie algebra structure defined by the tensor
$W[\alpha]_i^{jk}$ is equivalent to the direct sum of $m$ copies of
${\cal G}$ with canonical structure on it plus $n-m$ copies of ${\cal
G}$ with the structure of Abelian algebra over it.
\end{Proposition}

Let us also mention, see \cite{Ya3}, that the bracket $[X,Y]_A$
corresponds to the trivial cocycle of $\gl(n,{\cal A})$, - the Lie
algebra defined by the commutator on $\Mat({n,\cal A})$:
\begin{equation}\label{eq:51}
[X,Y]_A=d\beta_A(X,Y)=[X,\beta_A(Y)]-[Y,\beta_A(X)]-\beta_A([X,Y])
\end{equation}
where $\beta_A :{\cal G}^n\mapsto {\cal G}^n$ is given by
\begin{equation}\label{eq:52}
\beta_A(X)={1\over 2}(AX+XA).
\end{equation}
However, as we shall see below, the above facts does not depreciate the
construction we are going to present because the bracket $[{\bf
x},{\bf y}]^{e_0}$ is just a starting point, in what follows we shall
deform it and restrict it to subalgebras.

\section{Deformations of the bracket $[{\bf x},{\bf y}]^{e_0}$}

At the beginning let us consider the following result: Let 
$$
{\cal H}={\cal H}^{(0)}\oplus{\cal H}^{(1)}\oplus\ldots
\oplus {\cal H}^{(n-1)}
$$
be ${\bf Z}_n$-graded Lie algebra, that is:
\begin{equation}\label{eq:6}
[{\cal H}^{(i)},{\cal H}^{(j)}]\subset {\cal H}^{(i+j)(mod (n))}.
\end{equation}
If $x=\sum\limits_{i=1}^n x_i$, $y=\sum\limits_{i=1}^n y_i$, with
$x_i,y_i\in {\cal H}^{(i)}$ then
\begin{equation}\label{eq:7}
[x,y]_i=\sum\limits_{s+l=i~(mod(n))}[x_s,y_l].
\end{equation}
\begin{Proposition}\label{Proposition:2}
For any $\lambda\in {\bf K}$ the formula:
\begin{equation}\label{eq:8}
[x,y]^{\lambda}_i=
\sum\limits_{s+l=i}[x_s,y_l]+\lambda\sum\limits_{s+l\geq n}
[x_s,y_l], \qquad s+l=i (\mod n)
\end{equation}
defines a Lie algebra structure over ${\cal H}$.
\end{Proposition}
This proposition is a generalization of the following fact,
which is well known, and effectively used in order to construct
compatible Poisson brackets, see \cite{TrFo}.
\begin{Corollary}
Let ${\cal H}={\cal H}^{(0)}\oplus {\cal H}^{(1)}$ be a ${\bf Z}_2$
graded Lie algebra:
$$
[{\cal H}^{(0)},{\cal H}^{(0)}]\subset {\cal H}^{(0)},\quad
[{\cal H}^{(0)},{\cal H}^{(1)}]\subset {\cal H}^{(1)},\quad
[{\cal H}^{(1)},{\cal H}^{(1)}]\subset {\cal H}^{(0)}.
$$
Let $x=x_0+x_1$, $y=y_0+y_1$, $x_i,y_i\in {\cal H}^{(i)},~~i=0,1$.
Then for every $\lambda\in {\bf K}$ 
\begin{equation}\label{eq:9}
[x,y]^{\lambda}=[x_0,y_0]+\lambda[x_1,y_1]+ [x_0,y_1]+[x_1,y_0]
\end{equation}
is a Lie bracket.
\end{Corollary}
{\it Proof of the proposition.} Let us consider the commutative
algebra ${\cal B}$ over ${\bf K}$ generated by the elements $g^p,~
p=0,1,2,\ldots n-1$ with $g^n=0$. (We can take for example the
algebra of all the matrices of the type $\mu g^k$, $\mu\in {\bf K}$,
$k=0,1,2\ldots$ and $g$ the matrix with components
$g_{i}^{j}=\delta_{i+1}^j$). Then ${\cal B}$ is $n$-dimensional
algebra and ${\cal H}\otimes_{\bf K} {\cal B}$ has a natural
structure of Lie algebra. Let us construct the imbedding $h: {\cal
H}\mapsto {\cal H}\otimes_{\bf K} {\cal B}$ (as a vector space) 
defined as follows. For $x=\sum_{i=0}^{n-1} x_i$, $x_i\in {\cal H}^(i)$
\begin{equation}
h(x)=\sum\limits_{i=0}^{n-1} x_i\otimes
g^i.
\end{equation}
It is easy to see that the set $g({\cal H})$ is subalgebra in 
${\cal H}\otimes_{\bf K} {\cal B}$.
The bracket in it has the form: if $x=\sum_{i=1}^{n-1} x_i\otimes g^i
$, $y=\sum_{i=1}^{n-1} y_i\otimes g^i$ with $x_i,y_i\in {\cal H}^{(i)}$
then
\begin{equation}\label{eq:10}
[h(x),h(y)]=\sum\limits_{k=0}^{n-1}\left(\sum\limits_{i+j=k}[x_i,y_j]\right)
\otimes g^{k}.
\end{equation}
As $h({\cal H})$ is isomorphic to ${\cal H}$ the bracket
(\ref{eq:10}) induces a bracket on ${\cal H}$:
$[x,y]^0=h^{-1}([h(x),h(y)])$ for which $h$ is Lie-algebra isomorphism.
The resulting Lie-algebra bracket corresponds to the case $\lambda=0$ in 
(\ref{eq:8}). Therefore, (\ref{eq:8}) can be written in the following way:
\begin{equation}\label{eq:11}
\begin{array}{c}
[x,y]^{\lambda}=[x,y]^0+\lambda [x,y]^{a},\\[4pt]
([x,y]^{a})_i\equiv\sum\limits_{s+l\geq n}[x_s,y_l],\quad s+l=i~(\mod(n))
\end{array}
\end{equation}
and we have that for $\lambda=0$ and $\lambda=1$ the expression
$[x,y]^{\lambda}$ is a Lie bracket. From the other side $[x,y]^a$ is
Lie bracket too, because it is skewsymmetric and
\begin{equation}
\begin{array}{c}
([[x,y]^a,z]^a)_s+cycl(x,y,z)=\sum\limits_{i+j+k\geq
2n}([[x_i,y_j],z_k]+cycl(x,y,z))=0,\\[4pt]
i+j+k=s~(\mod(n))
\end{array}
\end{equation}
as a result of the Jacobi identity for ${\cal G}$. (Here the notation
$cycl(x,y,z)$ means that we must add to the corresponding 
expression two more terms, obtained from the first one by cyclic
permutation of $x,y,z$.

Finally, the proposition follows from the fact that $[x,y]^I, [x,y]^{II}$
are compatible Lie brackets. Q.E.D.

Consider now the algebra corresponding to the bracket $[{\bf x},{\bf
y}]^{e_0}$ over ${\cal G}^n$. This algebra has natural grading:
\begin{equation}
{\cal G}^n={\cal F}_n^{(0)}\oplus{\cal F}_n^{(1)}\oplus\ldots
\oplus {\cal F}_n^{(n-1)}
\end{equation}
where ${\cal F}_n^{(i)}$ consists of those ${\bf x}$ for which
$x_j=0$ for all $j\neq i$. We now use the proposition
(\ref{Proposition:2}) to modify the bracket $[{\bf x},{\bf
y}]^{e_0}$. We consider it as the bracket corresponding to
$\lambda=1$ and we take that corresponding to $\lambda=0$. In terms
of the tensor $W_{s}^{ij}$  we have $W_s^{ij}=\delta_{s}^{i+j}$,
$i,j,s=0,1,\ldots,n-1$. Let us denote the space ${\cal G}^n$ endowed
with our new algebra structure by ${\cal G}_L^{n}$. For it
\begin{equation}
\begin{array}{c}
[{\cal F}^{(i)}_n,{\cal F}^{(j)}_n]_L\subset {\cal F}_n^{(i+j)},\quad
i+j\leq n-1\\[4pt]
[{\cal F}^{(i)}_n,{\cal F}^{(j)}_n]_L \subset 0,\quad
i+j >n-1
\end{array}
\end{equation}
The extension we have obtained corresponds to the so called splitting
Leibnitz extension bracket, see \cite{JLT}. Its solvable part is
exactly the solvable Leibnitz extension:
$\bar{W}_{k}^{ij}=\delta^{i+j}_k$, $1\leq i,j,k\leq n-1$.

\section{Conclusion}
We have shown that some of the universal Lie algebra extensions can
be naturally understand within the frames of the theory of the
bundles of Lie brackets and their deformations. Whether we can obtain
all of them in this way is a very interesting open question.

\end{document}